# A matrix method based on the Fibonacci polynomials to the generalized pantograph equations with functional arguments

Ayşe Betül Koç*,a, Musa Çakmak b, Aydın Kurnaz a

\* Correspondence: aysebetulkoc @ selcuk.edu.tr
aDepartment of Mathematics, Faculty of Science, Selcuk University, Konya, Turkey
bMustafa Kemal University, Yayladagı Vocational School, Hatay, Turkey

**Abstract:** In this study, a collocation method based on the Fibonacci operational matrix is proposed to solve generalized pantograph equations with linear functional arguments. Some illustrative examples are given to verify the efficiency and effectiveness of the proposed method.

**Keywords:** Pantograph equations, Matrix method, Polynomial approximation, Operational matrices, Fibonacci polynomials

**MSC:** 65D25.

## 1. Introduction

Many phenomena in applied branches that fail to be modeled by the ordinary differential equations can be described by the Delay differential equations. Many researchers have studied different applications of those equations in variety of applied sciences such as biology, physics, economy and electrodynamics (see [1-4]). Pantograph equations with proportional delays play an important role in this context. The existence and uniqueness of the analytic solutions of the multi-pantograph equation are investigated in [5]. A numerical approach to multi-pantograph equations with variable coefficients is also studied in [6]. An extension of the multi-pantograph equation is known to be the generalized pantograph equation with functional arguments defined as

$$y^{(m)}(x) = \sum_{j=0}^{J}\sum_{k=0}^{m-1} p_{jk}(x) y^{(k)}(\alpha_{jk} x + \beta_{jk}) + g(x), \quad x \in I = [a,b] \tag{1.1}$$

under the mixed conditions

$$\sum_{k=0}^{m-1} c_{rk} y^{(k)}(\mu_{rk}) = \lambda_r, \quad r = 0(1)(m-1), \quad a \le \mu_{rk} \le b \tag{1.2}$$

where $\alpha_{jk}$, $\beta_{jk}$, $c_{rk}$, $\mu_{rk}$ and $\lambda_r$ are real and/or complex coefficients; $p_{jk}$ and $g(x)$ are the analytical functions defined in interval $a \le x \le b$.

In recent years, many researchers have developed different numerical approaches to the generalized pantograph equations as variational iteration method [7], differential transform approach [8], Taylor method [9], collocation method based on Bernoulli matrix [10] and Bessel collocation method [11]. In this study, we investigate a collocation method based on the Fibonacci polynomial operational matrix for the numerical solution of the generalized pantograph equation (1.1). Even the Fibonacci numbers have been known for a long time, the Fibonacci polynomials are very recently defined to be an important agent in the world of polynomials [12-13]. Compared to the methods of the orthogonal polynomials, the Fibonacci approach has proved to give more precise and reliable results in the solution of differential equations [14].

This study is organized as follows: In the second part, a short review of the Fibonacci polynomials is presented. A Fibonacci operational matrix for the solution of the pantograph equation is developed in Section 3. Some numerical examples are given in Section 4 to illustrate efficiency and effectiveness of the method.

## 2. Operational matrices of the Fibonacci polynomials

The Fibonacci polynomials $F_r(x)$ are determined by following general formula [12-13]

$$F_{r+1}(x) = xF_r(x) + F_{r-1}(x), \quad \text{for } r > 1, \tag{2.1}$$

with $F_1(x) = 1$ and $F_2(x) = x$. Now, we will mention some matrix relations in terms of Fibonacci polynomials.

### 2.1 Fibonacci Series Expansions

To obtain an expansion form of the analytic solution of the pantograph equation, we use the Fibonacci collocation method as follows:

Suppose that the equation (1.1) has a continuous function solution that can be expressed in the Fibonacci polynomials

$$y(x) = \sum_{r=1}^{\infty} a_r F_r(x). \tag{2.2}$$

Then, a truncated expansion of $N$-Fibonacci polynomials can be written in the vector form

$$y(x) \cong \sum_{r=1}^{N} a_r F_r(x) = \mathbf{F}(x) \mathbf{A} \tag{2.3}$$

where the Fibonacci row vector $\mathbf{F}(x)$ and the unknown Fibonacci coefficients column vector $\mathbf{A}$ are given, respectively, by

$$\mathbf{F}(x) = \begin{bmatrix} F_1(x) & F_2(x) & \cdots & F_N(x) \end{bmatrix}_{1 \times N} \quad (2.4)$$

and

$$\mathbf{A} = \begin{bmatrix} a_1 & a_2 & \cdots & a_N \end{bmatrix}^T_{N \times 1}. \quad (2.5)$$

### 2.2 Matrix relations of the derivatives:

The $k$-th order derivative of (2.3) can be written as

$$y^{(k)}(x) \cong \sum_{r=1}^{N} a_r^{(k)} F_r(x) = \mathbf{F}(x)\mathbf{A}^{(k)}, \quad k = 0,1,\ldots,n \quad (2.6)$$

where $a_r^{(0)} = a_r$, $y^{(0)}(x) = y(x)$ and

$$\mathbf{A}^{(k)} = \begin{bmatrix} a_1^{(k)} & a_2^{(k)} & \cdots & a_N^{(k)} \end{bmatrix}^T \quad (2.7)$$

is the coefficient vector of the polynomial approximation of $k$-th order derivative. Then, there exists a relation between the Fibonacci coefficients as

$$\mathbf{A}^{(k+1)} = \mathbf{D}^k \mathbf{A}, \quad k = 0,1,\ldots,n \quad (2.8)$$

where $\mathbf{D}$ is $N \times N$ operational matrix for the derivative defined by [14]

$$\mathbf{D} = [d_{i,j}] = \begin{cases} i \cdot \sin \dfrac{(j-i)\pi}{2}, & j > i \\ 0, & j \leq i \end{cases}. \quad (2.9)$$

Making use of Eq.s (2.6) and (2.8) yields

$$y^{(k)}(x) = \mathbf{F}(x)\mathbf{D}^k \mathbf{A}, \quad k = 0,1,\ldots,n. \quad (2.10)$$

### 3. Solution procedure for the pantograph differential equations

Let us recall the $m$-th order linear pantograph-differential equation,

$$y^{(m)}(x) = \sum_{j=0}^{J} \sum_{k=0}^{m-1} p_{jk}(x) y^{(k)}(\alpha_{jk} x + \beta_{jk}) + g(x), \quad x \in I = [a,b]. \quad (3.1)$$

The first step in the solution procedure is to define the collocation points in the domain $I$, so that,

$$x_i = a + \frac{b-a}{N-1}(i-1), \quad i = 1,2,\ldots,N, \quad a \leq x_i \leq b. \quad (3.2)$$

Then, collocating problem (3.1) at the points in (3.2) yields

$$y^{(m.)}(x_i) = \sum_{j=0}^{J}\sum_{k=0}^{m-1} P_{jk}(x_i) y^{(k)}(\alpha_{jk} x_i + \beta_{jk}) + g(x_i), \quad i = 1(1)N, \; m+1 \leq N \quad (3.3)$$

The system (3.3) can, alternatively, be rewritten in the matrix form

$$\mathbf{Y}^{(m)} - \sum_{j=0}^{J}\sum_{k=0}^{m-1} \mathbf{P}_{jk} \bar{\mathbf{Y}}_{jk}^{(k)} = \mathbf{G} \quad (3.4)$$

where

$$\mathbf{P}_{jk} = \begin{bmatrix} P_{jk}(x_1) & 0 & \cdots & 0 \\ 0 & P_{jk}(x_2) & \cdots & 0 \\ \vdots & \vdots & \ddots & \vdots \\ 0 & 0 & \cdots & P_{jk}(x_N) \end{bmatrix} \text{ and } \mathbf{G} = [g(x_1) \; g(x_2) \; \cdots g(x_N)]^T.$$

Therefore, the $k$-th order derivative of the unknown function at the collocation points can be written in the matrix form as

$$[y^{(k)}(x_i)] = \mathbf{F}(x_i) \mathbf{D}^k \mathbf{A}, \quad i = 1, 2, \ldots, N$$

or equivalently,

$$Y^{(k)} = \begin{bmatrix} y^{(k)}(x_1) \\ y^{(k)}(x_2) \\ \vdots \\ y^{(k)}(x_N) \end{bmatrix} = \mathbf{F} \mathbf{D}^k \mathbf{A}. \quad (3.5)$$

To express the functional terms of Eq. (1.1) as in the form (2.3), let $x \to \alpha_{jk} x_i + \beta_{jk}$ in the relation (3.5), then obtain

$$[y^{(k)}(\alpha_{jk} x_i + \beta_{jk})] = \bar{\mathbf{F}}(\alpha_{jk} x_i + \beta_{jk}) \mathbf{D}^k \mathbf{A}, \quad i = 1(1)N \quad (3.6)$$

or,

$$\bar{\mathbf{Y}}_{jk}^{(k)} = \begin{bmatrix} y^{(k)}(\alpha_{jk} x_1 + \beta_{jk}) \\ y^{(k)}(\alpha_{jk} x_2 + \beta_{jk}) \\ \vdots \\ y^{(k)}(\alpha_{jk} x_N + \beta_{jk}) \end{bmatrix} = \bar{\mathbf{F}}_{jk} \mathbf{D}^k \mathbf{A} \quad (3.7)$$

where $\bar{\mathbf{F}}_{jk}$ are Fibonacci operational matrices corresponding to the coefficients $\alpha_{jk}$. Then, replacing (3.5) and (3.7) in equation (3.4) gives the fundamental matrix equation for the problem (3.1) as

$$\left\{ \mathbf{F}\mathbf{D}^m - \sum_{j=0}^{J}\sum_{k=0}^{m-1} \mathbf{P}_{jk}\bar{\mathbf{F}}_{jk}\mathbf{D}^k \right\} \mathbf{A} = \mathbf{G} \tag{3.8}$$

which corresponds to a system of $(N)$ algebraic equations for the unknown Fibonacci coefficients $a_r$, $r = 1, 2, ..., N$. In other words, when we denote the expression in the sum by $\mathbf{W} = [w_{s,t}]$, for $s = 1, 2, ..., N$ and $t = 1, 2, ..., N$, we get

$$\mathbf{W}\mathbf{A} = \mathbf{G}. \tag{3.9}$$

Thus, the augmented matrix of Eq. (3.9) becomes

$$[\mathbf{W}; \mathbf{G}]. \tag{3.10}$$

On the other hand, with the help of Eq.(2.10), the conditions (1.2) can be converted to following matrix form

$$\mathbf{U}_j = \left[ \sum_{k=0}^{m-1} c_{rk} y^{(k)}(\mu_{rk}) \right] = [\lambda_r,] \tag{3.11}$$

where,

and

$$\mathbf{U}_j = [u_{\rho,\sigma}] = \sum_{k=0}^{m-1} c_{rk} \mathbf{F}_{\mu_{rk}} \mathbf{D}^k \mathbf{A},$$
$$\mathbf{F}_{\mu_{rk}} = [F_1(\mu_{rk}) \quad F_2(\mu_{rk}) \quad \cdots \quad F_N(\mu_{rk})]. \tag{3.12}$$

Therefore, the augmented matrix of the specified conditions is

$$[\mathbf{U}_j; \lambda_j] = [u_{j_1} \quad u_{j_2} \quad \cdots \quad u_{j_N} : \gamma_j]. \tag{3.13}$$

Consequently, (3.10) together with (3.13) can be written in the new augmented matrix form

$$[\mathbf{W}^* : \mathbf{G}^*]. \tag{3.14}$$

This form can also be achieved by replacing some rows of the matrix (3.10) by the rows of (3.13) or adding those rows to the matrix (3.10) provided that $\det(\mathbf{W}^*) \neq 0$. Finally, the vector $A$ (thereby vector of the coefficients $a_r$) is determined by applying some numerical methods designed especially to solve the system of linear equations. On the other hand, when the singular case $\det(\mathbf{W}^*) = 0$ appears, the least square methods are inevitably available to reach the best possible approximation. Therefore, the approximated solution can be obtained.

This would be the Fibonacci series expansion of the solution to the problem (3.1) with specified conditions.

**Accuracy of the Results:**

We can, now, proceed with a short accuracy analysis of the problem in a similar way to [15]. As the truncated Fibonacci series expansion is an approximate solution of Eq. (1.1) with (1.2), it must be satisfy the following equality for $x = x_r \in [a,b]$

$$E(x_r) = \left| y^{(m.)}(x_r) - \sum_{j=0}^{J}\sum_{k=0}^{m-1} P_{jk}(x_r) y^{(k)}(\alpha_{jk} x_r + \beta_{jk}) - g(x_r) \right| \cong 0$$

or

$$E(x_r) \leq 10^{-k_r} \quad (k_r \text{ is any positive integer}).$$

When $\max(10^{-k_r}) = 10^{-k}$ ($k$ is any integer) is prescribed, the truncation limit $N$ is increased until the difference $E(x_r)$ at each of the collocation points becomes smaller than the desired value $10^{-k}$.

**4. Numerical results**

In this part, three illustrative examples are given in order to clarify the findings of the previous section. The errors of the proposed method are compared with those of the errors occurred in the solutions by some other methods in Table 1-3 for two sample examples. It is noted here that the number of collocation points in the examples is indicated by the capital letter $N$.

**Example 1** [10-11] Consider the following linear pantograph equation

$$y''(x) = \frac{3}{4} y(x) + y\left(\frac{x}{2}\right) - x^2 + 2, \quad 0 \leq x \leq 1$$

with the initial conditions

$$y(0) = y'(0) = 0 .$$

The exact solution of this problem is known $y(x) = x^2$. When the solution procedure in Section 3 is applied to the problem, the solution of the linear algebraic system gives the numerical approximation of the solution to the problem. It is noteworthy that the method reaches the exact solution $y(x) = x^2$ even for $N = 3$.

**Example 2** [9,16-17] Now, consider the following equation with variable coefficient given in

$$\begin{cases} y'(x) = \frac{1}{2} e^{x/2} y\left(\frac{x}{2}\right) + \frac{1}{2} y(x), \ 0 \leq x \leq 1 \\ y(0) = 1. \end{cases}$$

The exact solution is also known to be $y(x) = \exp(x)$. A comparison of the absolute errors of the proposed approach, Taylor method [16] and the exponential approach [17] is given in Table1 for $N = 5$. Another comparison of the present method with the methods of Taylor polynomials [9,16] is also given for $N = 9$ and $N = 8$ in Table 2. These results verify that the Fibonacci approach is better at least one decimal place in accuracy than the others.

Table 1 Comparison of the absolute errors of different approximation techniques to Example 2

| $x$ | Present Method N=5 | Exponential approach [17] N=5 | Taylor polynomial approach [16] N=5 |
|---|---|---|---|
| 0.2 | 0.2553 E-05 | 0.32778 E-02 | 0.271 E-06 |
| 0.4 | 0.1965 E-05 | 0.32081 E-02 | 0.882 E-05 |
| 0.6 | 0.3874 E-05 | 0.44444 E-02 | 0.682 E-04 |
| 0.8 | 0.4833 E-05 | 0.46898 E-02 | 0.293 E-03 |
| 1.0 | 0.2690 E-04 | 0.12864 E-01 | 0.912 E-03 |

Table 2 Comparison of the absolute errors of different approximation techniques to Example

| $x$ | Present Method N=9 | Taylor polynomial approach [16] N=9 | Taylor Method [9] N=8 |
|---|---|---|---|
| 0.2 | 0 | 0.271 E-06 | 1.44 E-12 |
| 0.4 | 0 | 0.882 E-05 | 7.524 E-10 |
| 0.6 | 0 | 0.682 E-04 | 2.953 E-08 |
| 0.8 | 0.1 E-08 | 0.293 E-03 | 4.018 E-07 |
| 1.0 | 0.1 E-08 | 0.912 E-03 | 3.059 E-06 |

**Example 3** [6,18] Finally, let us consider the pantograph equation with variable coefficients

$$\begin{cases} y'(x) = -y(x) - e^{-0.5x} \sin(0.5x) y(0.5x) - 2e^{-0.75x} \cos(0.5x) \sin(0.25x) y(0.25x), \ 0 \leq x \leq 1 \\ y(0) = 1 \end{cases}$$

which has the exact solution $y(x) = e^{-x}\cos x$. Computed results are compared with the results of the Taylor [6] and Boubaker [18] matrix methods in Table 3.

Table 3 Comparison of the absolute errors of different approximation techniques to Example 3

| $x$ | Present Method | | | Taylor Matrix Method [6] | | Baubaker Matrix Method [18] | |
|---|---|---|---|---|---|---|---|
| | N=5 | N=9 | N=12 | N=5 | N=9 | N=9 | N=12 |
| 0.2 | 0.18903 E-5 | 0.12102 E-10 | 0.18730 E-14 | 0.69082 E-6 | 0.1300 E-8 | 0.121 E-10 | 0 |
| 0.4 | 0.62395 E-6 | 0.96855 E-11 | 0.14897 E-14 | 0.42924 E-4 | 0.1434 E-6 | 0.968 E-11 | 0.310 E-14 |
| 0.6 | 0.13542 E-5 | 0.71954 E-11 | 0.11151 E-14 | 0.47443 E-3 | 0.2058 E-5 | 0.719 E-11 | 0.110 E-14 |
| 0.8 | 0.15097 E-5 | 0.68229 E-11 | 0.68964 E-15 | 0.25855 E-2 | 0.1212 E-4 | 0.682 E-11 | 0.200 E-14 |
| 1.0 | 0.47735 E-4 | 0.75830 E-9 | 0.13256 E-12 | 0.95631 E-2 | 0.4003 E-4 | 0.758 E-9 | 0.562 E-12 |

**Acknowledgements**

This study was supported by Research Projects Center (BAP) of Selcuk University. Also, the authors would like to thank the Selcuk University and TUBITAK for their supports.

We have denoted here that a minor part of this study was presented orally at the "2nd International Eurasian Conference on Mathematical Sciences and Applications (IECMSA-2013)", Sarajevo, August, 2013.